\documentclass[12pt]{amsart}
\usepackage{amscd,amsmath,amsthm,amssymb,enumerate}
\usepackage[left]{lineno}
\usepackage{pstricks}
\usepackage[utf8]{inputenc}
\usepackage{epstopdf}

\usepackage{comment}

\usepackage[all]{xy}

 \def\Soc{{\mathbf Soc}}

 \def\opn#1#2{\def#1{\operatorname{#2}}} 
 \opn\chara{char} \opn\length{\ell} \opn\pd{pd} \opn\rk{rk}
 \opn\projdim{proj\,dim} \opn\injdim{inj\,dim} \opn\rank{rank}
 \opn\depth{depth} \opn\grade{grade} \opn\height{height}
 \opn\bigheight{bigheight}
 \opn\embdim{emb\,dim} \opn\codim{codim}
 
 \opn\superheight{superheight}\opn\lcm{lcm}
 \opn\trdeg{tr\,deg}
 \opn\reg{reg} \opn\lreg{lreg} \opn\ini{in} \opn\lpd{lpd}
 \opn\size{size} \opn\sdepth{sdepth}
 \opn\link{link}\opn\fdepth{fdepth}\opn\lex{lex}
 \opn\type{type}
 \opn\gap{gap}
 \opn\arithdeg{arith-deg}
 \opn\Deg{Deg}
 \opn\sat{sat}
 \opn\mat{mat}
 \opn\Mat{Mat}
 \opn\div{div} \opn\Div{Div} \opn\cl{cl} \opn\Cl{Cl}
 \opn\Spec{Spec} \opn\Supp{Supp} \opn\supp{supp} \opn\Sing{Sing}
 \opn\Ass{Ass} \opn\Min{Min}\opn\Mon{Mon} \opn\Max{Max}
 \opn\Ann{Ann} \opn\Rad{Rad} \opn\Soc{Soc}
 \opn\Im{Im} \opn\Ker{Ker} \opn\Coker{Coker} \opn\Am{Am}
 \opn\Hom{Hom} \opn\Tor{Tor} \opn\Ext{Ext} \opn\End{End}
 \opn\Aut{Aut} \opn\id{id}
 
 \opn\nat{nat}
 \opn\pff{pf}
 \opn\Pf{Pf} \opn\GL{GL} \opn\SL{SL} \opn\mod{mod} \opn\ord{ord}
 \opn\Gin{Gin} \opn\Hilb{Hilb}\opn\sort{sort}
 \opn\PF{PF}\opn\Ap{Ap}
 \opn\mult{mult}
 \opn\bight{bight}
 \opn\aff{aff}
 \opn\relint{relint} \opn\st{st}
 \opn\lk{lk} \opn\cn{cn} \opn\core{core} \opn\vol{vol}  \opn\inp{inp} \opn\nilpot{nilpot}
 \opn\link{link} \opn\star{star}\opn\lex{lex}\opn\set{set}
 \opn\width{wd}
 \opn\Fr{F}
 \opn\QF{QF}
 \opn\G{G}
 \opn\type{type}\opn\res{res}
 \opn\conv{conv}
 \opn\Shad{Shad}
 \opn\gr{gr}
 
 \def\pot#1#2{#1[\kern-0.28ex[#2]\kern-0.28ex]}

 \opn\dirlim{\underrightarrow{\lim}}
 \opn\inivlim{\underleftarrow{\lim}}
 \let\to=\rightarrow
 
 \def\Implies{\ifmmode\Longrightarrow \else
         \unskip${}\Longrightarrow{}$\ignorespaces\fi}
 \def\implies{\ifmmode\Rightarrow \else
         \unskip${}\Rightarrow{}$\ignorespaces\fi}
 \def\iff{\ifmmode\Longleftrightarrow \else
         \unskip${}\Longleftrightarrow{}$\ignorespaces\fi}

 \let\:=\colon
\theoremstyle{plain}
\newtheorem{theorem}{Theorem}[section]
\newtheorem{Theorem}[theorem]{Theorem}

\newtheorem{Proposition}[theorem]{Proposition}

\newtheorem{Corollary}[theorem]{Corollary}

\newtheorem{Lemma}[theorem]{Lemma}
\newtheorem{lem}[theorem]{Lemma}
\newtheorem{claim}{Claim}

\theoremstyle{definition}

\newtheorem{Definition}[theorem]{Definition}

\newtheorem{Example}[theorem]{Example}

\newtheorem{Remark}[theorem]{Remark}

\newtheorem{Fact}[theorem]{Fact}
\newtheorem*{acknowledgments}{Acknowledgments}

 \let\epsilon\varepsilon
 \let\kappa=\varkappa
 \def\qed{\ifhmode\textqed\fi
       \ifmmode\ifinner\quad\qedsymbol\else\dispqed\fi\fi}
 \def\textqed{\unskip\nobreak\penalty50
        \hskip2em\hbox{}\nobreak\hfil\qedsymbol
        \parfillskip=0pt \finalhyphendemerits=0}
 \def\dispqed{\rlap{\qquad\qedsymbol}}
 
 \opn\dis{dis}
 \def\pnt{{\raise0.5mm\hbox{\large\bf.}}}
 
 \opn\Lex{Lex}

\newcommand{\rme}{\mathrm{e}}

\newcommand{\rmr}{\mathrm{r}}

\newcommand{\rmK}{\mathrm{K}}

\newcommand{\rmQ}{\mathrm{Q}}

\newcommand{\fka}{\mathfrak{a}}

\newcommand{\fkm}{\mathfrak{m}}
\newcommand{\fkn}{\mathfrak{n}}

\newcommand{\fkp}{\mathfrak{p}}

\newcommand{\fkM}{\mathfrak{M}}

\def\ol{\overline}
\def\canred{\mathrm{can. red\,}}
\def\tr{\mathrm{tr}}

\title{The reduction number of canonical ideals}
\author{Shinya Kumashiro}
\address{Shinya Kumashiro: National Institute of Technology, Oyama College
771 Nakakuki, Oyama, Tochigi, 323-0806, Japan}
\email{skumashiro@oyama-ct.ac.jp}

\thanks{2020 {\em Mathematics Subject Classification.} Primary: 13H10, Secondary: 13B02, 13D40}
\thanks{{\em Key words and phrases.} Cohen-Macaulay ring, Gorenstein ring, almost Gorenstein ring, nearly Gorenstein ring, canonical ideal}
\thanks{The author was supported by JSPS KAKENHI Grant Numbers JP19J10579 and JP21K13766.}

\begin{document}

\begin{abstract}
In this paper, we introduce an invariant of Cohen-Macaulay local rings in terms of the reduction number of canonical ideals. The invariant can be defined in arbitrary Cohen-Macaulay rings and it measures how close to being Gorenstein. First, we clarify the relation between almost Gorenstein rings and nearly Gorenstein rings by using the invariant in dimension one. We next characterize the idealization of trace ideals over Gorenstein rings in terms of the invariant. It provides better prospects for a result of the almost Gorenstein property of idealiztion.
\end{abstract}

\maketitle



\section{Introduction}\label{section1}

The aim of this paper is to develop a theory of non-Gorenstein Cohen-Macaulay rings. Gorenstein rings are an important class of Cohen-Macaulay rings and known to have good properties such as duality and reflexivity. On the other hand, among the huge class of Cohen-Macaulay rings, the class of Gorenstein rings seems narrower.
For instance, although any normal semigroup rings are Cohen-Macaulay, a normal semigroup ring is Gorenstein only the case that its interior coincides with itself after some shift (see \cite[Theorem 6.3.5]{BH}).
Moreover, if $R$ is a Cohen-Macaulay local ring and $M$ is a maximal Cohen-Macaulay $R$-module, the idealization $R\ltimes M$ becomes Cohen-Macaulay local ring again. However, $R\ltimes M$ becomes Gorenstein only in the case that $M$ is the canonical module of $R$ (\cite{R}).
Therefore, it seems natural to expect a new class of rings between Gorenstein and Cohen-Macaulay.

{\it Almost Gorenstein rings} are one of the most interesting objects in the study of non-Gorenstein rings. The basic papers \cite{BF, GMP, GTT} revealed the properties of the non-Gorenstein almost Gorenstein local rings such as G-regularity and the Gorensteinness of the blow-up at the maximal ideals in dimension one. 
Besides the almost Gorenstein theory, the study of non-Gorenstein Cohen-Macaulay rings has been carried out under intense competition. One can also find other stratifications of Cohen-Macaulay rings in \cite{CGKM, DS, DS2, GK, HHS}. Especially, {\it nearly Gorenstein rings} have gained attention in recent years along with the notion of the trace ideal of the canonical module (\cite{DKT, HHS, MS}). 

In this paper we introduce a new invariant, say {\it canonical reduction number}. This invariant is, roughly to say, the reduction number of a canonical ideal. More precisely, we define the invariant by the new notion, say {\it almost reduction}. The notion of almost reduction is a natural generalization of reduction, but it will be useful to avoid the assumption that the residue field is infinite. 

The canonical reduction number can be defined in arbitrary Cohen-Macaulay rings and it measures how close to being Gorenstein. 
In this paper, we establish two theorems on the canonical reduction number. The first theorem builds bridges between almost Gorenstein rings and nearly Gorenstein rings as follows. 

\begin{Theorem}{\rm (Theorem \ref{qqq3.11})}\label{zzz1.1}
Let $R$ be a one-dimensional Cohen-Macaulay local ring possessing the canonical module. Then $R$ is almost Gorenstein if and only if $R$ is nearly Gorenstein and the canonical reduction number of $R$ is less than or equals to two.
\end{Theorem}

The second theorem is about idealization. 
Let $R$ be a commutative ring and $M$ an $R$-module. Let $A=R \ltimes M$ denote the {\it idealization} of $M$ over $R$, i.e., $A = R \oplus M$ as an $R$-module and the multiplication in $A$ is given by
$$(a,x)(b,y) = (ab, bx + ay)$$
where $a,b \in R$ and $x,y \in M$. Then it is known that the properties of $R\ltimes M$ are characterized by the properties of $R$ and $M$. Especially, if $R$ is a Cohen-Macaulay ring and $M$ is a maximal Cohen-Macaulay $R$-module, then the idealization $R\ltimes M$ is a Cohen-Macaulay ring. However, $R\ltimes M$ is Gorenstein only in the case that $M$ is isomorphic to the canonical module (\cite{R}). With this perspective, the question of when $R\ltimes M$ is close to Gorenstein is studied (\cite{GK2, GK, GKL, GMP, GTT}).
In this paper, we characterize the idealization of trace ideals of Gorenstein rings by using the canonical reduction number, which is a generalization of \cite[Corollary 6.4]{GMP}.
Let us recall that, for an $R$-module $M$, the {\it trace ideal} of $M$ is the image of the $R$-linear map
\[
\Hom_R(M, R) \otimes_R M \to R, \quad f\otimes x \mapsto f(x)
\]
where $f\in \Hom_R(M, R)$ and $x\in M$. An ideal $I$ is called a {\it trace ideal} if $I$ is a trace ideal of some $R$-module. One can consult \cite{GIK2, Lin, Lin2}  for the basic properties of trace ideals.

\begin{Theorem}{\rm (Theorem \ref{d2.3})}\label{zzz1.2}
Let $R$ be a Gorenstein local ring of dimension $d>0$ and $M$ be a maximal Cohen-Macaulay faithful $R$-module. Let $A=R\ltimes M$ be the idealization of $M$ over $R$, and $\omega_A$ denotes the canonical module of $A$. Then the following assertions are equivalent:
\begin{enumerate}[{\rm (a)}] 
\item $M$ is isomorphic to some trace ideal of $R$;
\item $\Hom_R(M, R)$ is isomorphic to some finite birational extension $B$ of $R$ such that $B_\fkM$ is a Cohen-Macaulay local ring of dimension $d$ for all $\fkM\in \Max B$;
\item the canonical reduction number of $A$ is less than or equals to two.
\end{enumerate}
When this is the case, if $A$ is not Gorenstein, then $\rmr(A)=\rmr(R/I)+2$ where $I$ denotes the trace ideal isomorphic to $M$.
\end{Theorem}

It is known that, for a one-dimensional Gorenstein local ring $(R, \fkm)$ and a maximal Cohen-Macaulay faithful $R$-module $M$, $R\ltimes M$ is almost Gorenstein if and only if either $M\cong R$ or $M\cong \fkm$ (\cite[Corollary 6.4]{GMP}). Since $R$ and $\fkm$ are trace ideals if $R$ is not a discrete valuation ring, Theorem \ref{zzz1.2} provides better prospects of \cite[Corollary 6.4]{GMP}.

Let us explain how constructed this paper. In Section \ref{section2} we first note the notion of almost reduction. This notion is a simple generalization of reduction, but it is effective to avoid the condition that the residue field is infinite. We then define the canonical reduction number. After that, we investigate the basic properties of the canonical reduction number, and characterize rings whose canonical reduction number is less than or equals to two.
In Section \ref{newnewsection2} we focus on the case of dimension one. In dimension one, we can define the Hilbert function of canonical ideals and we give a characterization of the canonical reduction number in terms of the Hilbert function. Theorem \ref{zzz1.1} is also shown in this section.
The purpose of Section \ref{newsection2}  is to show Theorem \ref{zzz1.2}.


Let us fix our notation. Throughout this paper, all rings are commutative Noetherian rings with identity. For a  ring $R$, $\rmQ(R)$ (resp. $\ol{R}$) denotes the total ring of fraction of $R$ (resp. the integral closure of $R$). For an $R$-module $M$, $M^*$  denotes the $R$-dual $\Hom_R (M, R)$. $\ell_R (M)$ denotes the length of $M$.

Suppose that $(R, \fkm)$ is a Noetherian local ring and $I$ is an $\fkm$-primary ideal of $R$. Then $\ell_R (R/I^n)$ agrees with a polynomial function of degree $d=\dim R$ for all $n\gg 0$. We then write
{\small
\[
\ell_R (R/I^{n+1})=\rme_0(I) \binom{n+d}{d} -\rme_1(I) \binom{n+d-1}{d-1} + \cdots +(-1)^{d-1} \rme_{d-1}(I) \binom{n+1}{1}  +(-1)^{d} \rme_{d}(I)
\]
}
with some integers $\rme_0(I)$, $\dots$, $\rme_d(I)$. The integers $\rme_0(I)$, $\dots$, $\rme_d(I)$ are called the {\it Hilbert coefficients} of $I$. If $R$ is a Cohen-Macaulay local ring, $\rmr(R)$ denotes the Cohen-Macaulay type of $R$. 

We say that an $R$-module $I$ is a {\it fractional ideal}, if $I$ is a finitely generated $R$-submodule of $\rmQ(R)$ containing a non-zerodivisor of $R$. For fractional ideals $I$ and $J$, $I:J$ (resp. $I:_R J$) stands for the set
\[
\left\{ \alpha \in \rmQ(R) \mid \alpha J\subseteq I \right\}
\]
(resp. $(I:J)\cap R=\left\{ \alpha \in R \mid \alpha J\subseteq I \right\}$). We freely use the following facts.

\begin{Remark}(\cite{HK})\label{b1.1}
Let $R$ be a Noetherian ring. Let $I$ and $J$ be fractional ideals. Then we have the following.
\begin{enumerate}[{\rm (a)}] 
\item If $I\cong J$, then $J=\alpha I$ for some $\alpha \in \rmQ(R)$. 
\item $I:J\cong \Hom_R (J, I)$, where $\alpha\in I:J$ corresponds to the multiplication map by $\alpha$.
\end{enumerate}
\end{Remark}


\section{almost reduction and canonical reduction number}\label{section2}

Let $R$ be a Noetherian ring. First of all, we note the notion of almost reduction, which is a simple extension of reduction. 

\begin{Definition}\label{f1.1}
Let $I$ and $J$ be  ideals of $R$. Then we say that $J$ is an {\it almost reduction} of $I$ if there exists an integer $\ell\ge 0$ such that
\begin{enumerate}[{\rm (a)}] 
\item $J^\ell \subseteq I^\ell$ and
\item $I^{\ell+1} = JI^{\ell}$.
\end{enumerate}
\end{Definition}


\begin{Remark}\label{h1.2}
Let $I$ and $J$ be ideals of $R$. Suppose that $J$ is an almost reduction of $I$. Then we have the following.
\begin{enumerate}[{\rm (a)}] 
\item $\sqrt{I}=\sqrt{J}$, where $\sqrt{I}$ denotes the radical of $I$. Hence $\height_R I=\height_R J$.
\item If $J^\ell \subseteq I^\ell$ and  $I^{\ell+1} = JI^{\ell}$ for $\ell\ge 0$, then 
$J^n \subseteq I^n$ and  $I^{n+1} = JI^{n}$ for all $n\ge \ell$.
\item Let $(R, \fkm)$ be a Noetherian local ring. If $I$ and $J$ are $\fkm$-primary ideals of $R$, then $\rme_0(J)=\rme_0(I)$.
\end{enumerate}
\end{Remark}

It is a well-known fact that, for a Noetherian local ring $(R, \fkm)$ with the infinite field $R/\fkm$, each $\fkm$-primary ideal $I$ has a parameter ideal $Q\subseteq I$ as its reduction. However, it is not true if $R/\fkm$ is finite. The following is an example of an $\fkm$-primary ideal which has no parameter reduction, but has a parameter almost reduction.

\begin{Example}{\rm (cf. \cite[Remark 2.10]{GMP})}\label{qwqw2.3}
Let $k[[X, Y, Z]]$ be the formal power series ring over the field $k=\mathbb{Z}/2\mathbb{Z}$. Set $R=k[[X, Y, Z]]/\fka$, where 
\[
\fka=(X, Y)\cap (Y, Z)\cap (Z, X)=(XY, YZ, ZX).
\]
Set $I=(x+y, y+z)$, where $x$, $y$, and $z$ denote the images of $X$, $Y$, and $Z$ in $R$ respectively. Then $(x+y+z)$ is an almost reduction of $I$, but $I$ has no parameter reduction.
\end{Example}

\begin{proof}
Set $\fkm=(x, y, z)$ and $f=x+y+z$. Then $I\ne \fkm=I+(f)$, whence $f\not\in I$. A standard calculation shows $I^n=(x^n, y^n, z^n)$ for all $n\ge 2$. Hence we have $fI^2=I^3$ and $f^2\in I^2$, which implies that $(f)$ is an almost reduction of $I$. Assume that $(a)$ is a reduction of $I$. Write $a=c_1 x + c_2 y + c_3 z + g$, where $c_1, c_2, c_3\in k=\mathbb{Z}/2\mathbb{Z}$ and $g\in \fkm^2$. Then we obtain 
\[
aI^n=(c_1x^{n+1}, c_2 y^{n+1}, c_3z^{n+1})=I^{n+1}
\]
for $n\gg 0$. Therefore, $c_1=c_2=c_3=1$. It follows that $f\in I$ since $a=f+g\in I$, thus it is a contradiction. 
\end{proof}

In Section \ref{newnewsection2} we explore the existence of almost reduction for $\fkm$-primary ideals in one-dimensional Cohen-Macaulay local rings. In what follows, let us focus on the case where the height of ideals is one. 

\begin{Proposition}\label{f1.2}
Let $R$ be a Noetherian ring. Let $I$ and $J$ be ideals of $R$ containing a non-zerodivisor of $R$. Suppose that $(a)$ and $(b)$ are almost reductions of $I$ and $J$ respectively. If $I\cong J$, then we have the following.
\begin{enumerate}[{\rm (a)}] 
\item For any $\ell\ge 0$, $(a)^\ell \subseteq I^\ell$ and $I^{\ell+1}=aI^\ell$ if and only if $R[\frac{I}{a}]=(\frac{I}{a})^\ell$, where $\frac{I}{a}=\{\frac{x}{a}\in \rmQ(R) \mid x\in I\}$ denotes a fractional ideal of $R$.
\item $R[\frac{I}{a}]=R[\frac{J}{b}]$ in $\rmQ(R)$.
\item For any $\ell\ge 0$, $(a)^\ell \subseteq I^\ell$ and $I^{\ell+1}=aI^\ell$ if and only if $(b)^\ell \subseteq J^\ell$ and $J^{\ell+1}=bJ^\ell$.
\end{enumerate}
\end{Proposition}

\begin{proof}
Note that $a$ is a non-zerodivisor of $R$ since $I^{n+1}=aI^n\subseteq (a)$ for $n\gg 0$.
Set $L_1=\frac{I}{a}$ and $L_2=\frac{J}{b}$ in $\rmQ(R)$. Then $R\subseteq L_1^n=L_1^{n+1}$ and $R\subseteq L_2^n=L_2^{n+1}$ for $n\gg0$. Furthermore, for $1\le t \le n-1$, we have $R{\cdot}L_1^t\subseteq L_1^{n+t}=L_1^{n}$. Hence $R[L_1]=L_1^n$ and $R[L_2]=L_2^n$.

On the other hand, by Remark \ref{b1.1} (1), we have $L_2=\alpha L_1$ for some $\alpha \in \rmQ(R)$. It follows that $R[L_1]=\alpha R[L_1]$ by substituting $L_2=\alpha L_1$ to $L_2^n=L_2^{n+1}$.

(a): The above argument actually shows that the ``only if'' part. Conversely, suppose that $R[L_1]=L_1^\ell$. If $\ell=0$, then $I\subseteq (a)$. Hence $I=aI_1$ for some ideal $I_1$. Then $(a)^n\subseteq I^n=a^n I_1^n\subseteq (a)^n$ for $n\gg 0$, whence $I_1=R$. Thus $I=(a)$.
Assume $\ell>0$. Then $L_1^{\ell+1}$ and $L_1^{\ell-1}$ are in $R[L_1]=L_1^\ell$. The latter implies that $L_1^\ell\subseteq L_1^{\ell+1}$, thus $L_1^\ell=L_1^{\ell+1}$. It follows that $I^{\ell+1}=aI^\ell$. Furthermore  we have $(a)^\ell \subseteq I^\ell$ since  $R\subseteq R[L_1]=L_1^\ell$.

(b): It follows from the observation that $R[L_2]=L_2^n=\alpha^n L_1^n=\alpha^n R[L_1]=R[L_1]$.

(c): We have only to show the ``only if'' part. By the assumption and (a) we have $L_1^\ell=R[L_1]$. Hence 
\[
L_2^\ell=\alpha^\ell L_1^\ell=\alpha^\ell R[L_1]=R[L_1]=R[L_2]
\]
by (b). It follows that $(b)^\ell \subseteq J^\ell$ and $J^{\ell+1}=bJ^\ell$ by (a).
\end{proof}

Proposition \ref{f1.2} (c) claims that the almost reduction number of $I$ and that of $J$ are equal. This fact provides an invariant of Cohen-Macaulay local rings.
Let us recall the following fact.

\begin{Fact} (\cite[Proposition 3.3.18]{BH})\label{f1.5}
Let $(R, \fkm)$ be a Cohen-Macaulay local ring possessing the canonical module $\omega_R$. Then the following conditions are equivalent:
\begin{enumerate}[{\rm (a)}] 
\item $R$ is generically Gorenstein, that is, $R_\fkp$ is Gorenstein for all $\fkp \in \Min R$;
\item $\omega_R$ has a rank;
\item there exists an ideal $\omega\subseteq R$ such that $\omega\cong \omega_R$.
\end{enumerate}
When this is the case, if $\omega\subsetneq R$, then $\height_R \omega=1$ and $R/\omega$ is Gorenstein.
\end{Fact}

We call an ideal $\omega\subseteq R$ {\it canonical} if $\omega$ is isomorphic to the canonical module of $R$. 

\begin{Definition}
Let $(R, \fkm)$ be a Cohen-Macaulay local ring possessing the canonical module. Suppose that $R$ is generically Gorenstein. We then call 
\[
\inf \left\{ n\ge 0 \ \middle| \ 
\begin{matrix}
\text{there exist a canonical ideal $\omega$ and an almost reduction $(a)$ of $\omega$}\\
\text{such that $\omega^{n+1}=a \omega^n$}
\end{matrix}\right\}
\]
the {\it canonical reduction number} of $R$, and denote by $\canred R$.
\end{Definition}


\begin{Remark}\label{lll1.7}
Let $(R, \fkm)$ be a Cohen-Macaulay local ring possessing the canonical module. Suppose that there exist a canonical ideal $\omega$ and an almost reduction $(a)$ of $\omega$. Then we have the following.
\begin{enumerate}[{\rm (a)}] 
\item For $n\ge 0$, $\omega^{n+1}=a\omega^n$ implies that $(a)^n\subseteq \omega^n$.
\item $\canred R_\fkp \le \canred R$ for all $\fkp \in \Spec R$.
\item Suppose that $x\in \fkm$ is a non-zerodivisor  of $R$ and $R/\omega$. Then $\canred R/xR\le \canred R$.
\end{enumerate}
\end{Remark}

\begin{proof}
(a) follows from $\left(\frac{\omega}{a}\right)^{n}=\left(\frac{\omega}{a}\right)^{m}=R[\frac{\omega}{a}]\supseteq R$ for $m\gg n$ by Proposition \ref{f1.2} (a). 
(b) follows from the fact that $\omega R_\fkp$ is a canonical ideal of $R_\fkp$. (c) follows from the fact that $(\omega +(x))/(x)$ is a canonical ideal of $R/(x)$ since $(\omega +(x))/(x)\cong \omega/((x) \cap \omega)=\omega/x\omega$.
\end{proof}

Flat local homomorphisms of rings also preserves the canonical reduction number.

\begin{Proposition}\label{i1.7}
Let $(R, \fkm) \to (S, \fkn)$ be a flat local homomorphism of Cohen-Macaulay rings such that $S/\fkm S$ is Gorenstein. Suppose that there exists the canonical module $\omega_R$ of $R$ and $R$ is generically Gorenstein. If $\canred R<\infty$, then $\canred R=\canred S$.
\end{Proposition}

\begin{proof}
Note that we have $\omega_S\cong S\otimes_R \omega_R$ since $S/\fkm S$ is Gorenstein. Hence, if $\omega$ is a canonical ideal of $R$, then $\omega S$ is a canonical ideal of $S$. 

Suppose that $\canred R<\infty$. Then there exist a canonical ideal $\omega$ and its almost reduction $(a)$. Let $n\ge 0$. Then we obtain that 
\begin{align*}
n\ge \canred R &\Leftrightarrow R\left[\textstyle\frac{\omega}{a}\right]/\left(\textstyle\frac{\omega}{a}\right)^n=0 &&\Leftrightarrow S\otimes_R \left(R\left[\textstyle\frac{\omega}{a}\right]/\left(\textstyle\frac{\omega}{a}\right)^n\right)=0\\
&\Leftrightarrow S\left[\textstyle\frac{\omega S}{a}\right]/\left(\textstyle\frac{\omega S}{a}\right)^n=0 && \Leftrightarrow n \ge \canred S
\end{align*}
by Proposition \ref{f1.2} (a). Therefore, we have $\canred R=\canred S$.
\end{proof}


The canonical reduction number measures how close to being Gorenstein.
In what follows, unless otherwise stated, let $(R, \fkm)$ be a Cohen-Macaulay local ring possessing the canonical module $\omega_R$. Set $d=\dim R$.

\begin{Proposition}\label{d1.1}
The following conditions are equivalent:
\begin{enumerate}[{\rm (a)}] 
\item $R$ is Gorenstein;
\item $\canred R= 0$;
\item $\canred R\le 1$;
\item there exist a canonical ideal $\omega \subseteq R$ and $a \in R$ such that $\omega^2=a\omega$.
\end{enumerate}
\end{Proposition}

\begin{proof}
Since (a) $\Rightarrow$ (b) $\Rightarrow$ (c) $\Rightarrow$ (d) is trivial, we have only to show that (d) $\Rightarrow$ (a). Note that $\omega^2=a\omega$ implies that 
$\frac{\omega}{a}\subseteq \frac{\omega}{a} : \frac{\omega}{a}=R$. 
Therefore, $\omega_1=\frac{\omega}{a}$ is a canonical ideal and $\omega_1^2=\omega_1$. 
Since $\omega_1$ is nonzero, it forces that $\omega_1=R$, whence $R$ is Gorenstein.
\end{proof}

\begin{Corollary}
Suppose that $R$ is a Cohen-Macaulay local normal domain. Then $\canred R<\infty$ if and only if $R$ is Gorenstein.
\end{Corollary}

\begin{proof}
We have only to show that the only if part. If $\canred R<\infty$, then we can choose a canonical ideal $\omega$ and its reduction $(a)$. By Proposition \ref{f1.2}, we have $R\subseteq R[\frac{\omega}{a}]\subseteq \ol{R}$. Thus $\canred R=0$.
\end{proof}

On the other hand, if $R$ is not a normal domain, then there are many non-Gorenstein Cohen-Macaulay rings such that the canonical redution number is finite, see Sections \ref{newnewsection2} and \ref{newsection2}.
We next characterize those rings $R$ with $\canred R\le 2$. To state our theorem, let us recall the definition of trace ideals.

\begin{Definition}
Let $R$ be a commutative ring. For an $R$-module $M$, the image of the evaluation map 
\begin{small}
\[
\varphi: \Hom_R(M, R)\otimes_R M \to R, \quad \text{where $\varphi(f\otimes x)=f(x)$ for $f\in \Hom_R(M, R)$ and $x\in M$},
\]
\end{small}
is called the {\it trace ideal} of $M$ and denoted by $\tr_R (M)$. We say that an ideal $I$ is a {\it trace ideal} of $R$ if $I=\tr_R (M)$ for some $R$-module $M$.
\end{Definition}

\begin{Remark}\label{b1.7}
\begin{enumerate}[{\rm (a)}] 
\item If $I$ is a fractional ideal, then $\tr_R(I)=(R:I)I$.
\item (\cite[Corollary 2.2]{GIK2}) Let $R$ be a Noetherian ring and $I$ an ideal containing a non-zerodivisor of $R$. Then the following conditions are equivalent:
\begin{enumerate}[{\rm (i)}] 
\item $I$ is a trace ideal of $R$;
\item $I=(R:I)I$, that is, $I=\tr_R (I)$;
\item $I:I=R:I$.
\end{enumerate}
\item (\cite[Lemma 2.1]{HHS}) The trace ideal $\tr_R (\omega_R)$ of the canonical module describes the non-Gorenstein locus of $R$, that is, 
\[
\left\{ \fkp\in \Spec R \mid \text{$R_\fkp$ is not Gorenstein} \right\}
=\left\{ \fkp\in \Spec R \mid \tr_R (\omega_R)\subseteq \fkp \right\}.
\]
\end{enumerate}
\end{Remark}





The following is a characterization of rings with canonical reduction number two.

\begin{Theorem}\label{d1.2}
The following conditions are equivalent:
\begin{enumerate}[{\rm (a)}] 
\item $\canred R \le 2$;
\item $\tr_R (\omega_R)\cong \omega_R^{*}$.
\end{enumerate}
When this is the case, letting $\omega$ be a canonical ideal of $R$ and $a\in \omega$ be an almost reduction of $\omega$, then $\tr_R(\omega_R)=R:R[\frac{\omega}{a}]$ and $R[\frac{\omega}{a}]=R: \tr_R(\omega_R)$.
\end{Theorem}

We prepare a lemma to prove Theorem \ref{d1.2}.

\begin{Lemma}\label{d1.5}
If $d=0$, then $\tr_R (\omega_R)\cong \omega_R^{*}$ if and only if $R$ is Gorenstein.
\end{Lemma}

\begin{proof}
If $R$ is Gorenstein, then $\tr_R (\omega_R)=R \cong \omega_R^{*}$ by Remark \ref{b1.7}. Conversely, suppose $\tr_R (\omega_R)\cong \omega_R^{*}$. Note that $\omega_R^{*}\cong \Hom_R (\omega_R, \Hom_R (\omega_R, \omega_R))\cong \Hom_R (\omega_R\otimes_R~\omega_R, \omega_R)$. Hence, by applying the $\omega_R$-dual to $\tr_R (\omega_R)\cong \omega_R^{*}$, we have 
\begin{align}\label{ddd}
\Hom_R (\tr_R (\omega_R), \omega_R) \cong \omega_R\otimes_R \omega_R.
\end{align}
On the other hand, by applying the $\omega_R$-dual to the exact sequence $0 \to \tr_R (\omega_R) \to R \to R/\tr_R (\omega_R) \to 0$, we have a surjection 
\begin{align}\label{ddd2}
\omega_R \to \Hom_R (\tr_R (\omega_R), \omega_R).
\end{align}
Therefore, from (\ref{ddd}) and (\ref{ddd2}), we obtain the surjection $\omega_R \to \omega_R \otimes_R \omega_R$. It follows that $\rmr(R)\ge \rmr(R)^2$, whence $R$ is Gorenstein.
\end{proof}

The following is a direct consequence of Lemma \ref{d1.5}.

\begin{Corollary}\label{d1.6}
If $\tr_R (\omega_R)\cong \omega_R^{*}$, then $R$ is generically Gorenstein. 
\end{Corollary}


\begin{proof}[Proof of Theorem \ref{d1.2}]
By Corollary \ref{d1.6}, we may assume that $R$ is generically Gorenstein. Since the conditions (a) and (b) are satisfied if $R$ is Gorenstein by Proposition \ref{d1.1} and Remark \ref{b1.7}(c), we may assume that $d>0$.

(a) $\Rightarrow$ (b): Choose a canonical ideal $\omega\subsetneq R$ and $a\in R$ so that $\omega^3=a\omega^2$ and $a^2\in \omega^2$.
Set $K=\frac{\omega}{a}$. Then we have  $R[K]=K^2$ by Proposition \ref{f1.2} (a). We obtain that
\begin{align*}
R:R[K]=R:K^2=(K:K):K^2=K:K^3=K:K^2=(K:K):K=R:K.
\end{align*}
It follows that 
\begin{align*}
(R:K)K=(R:R[K])K\subseteq (R:R[K])R[K] = R:R[K]=R:K.
\end{align*}
It follows that $(R:K)K^2\subseteq (R:K)K$, whence $R:K\subseteq (R:K)K$ since $(R:K)K^2=(R:R[K])R[K]=R:R[K]=R:K$. Therefore, we get $\tr_R (\omega_R)=(R:K)K=R:~K\cong \omega_R^*$ by Remark \ref{b1.7}(a).

(b) $\Rightarrow$ (a): Due to Corollary \ref{d1.6}, we may assume that there exists a canonical ideal $\omega\subsetneq R$. Then we have $\tr_R(\omega_R)=(R:\omega)\omega$ and $\omega_R^*\cong R:\omega$. Hence there exists $a\in \rmQ(R)$ such that $(R:\omega)\omega=a (R:\omega)$. We can replace $\omega$ and $a$ so that $a\in R$. Note that $(R:\omega)\omega^n=a^n (R:\omega)$ for all $n>0$.
We then obtain that 
\begin{align*}
\left(\textstyle\frac{\omega}{a}\right)^n&\subseteq(R:\omega):(R:\omega)=R:(R:\omega)\omega=(\omega:\omega):a(R:\omega)\\
&=\omega:a(R:\omega)\omega=\omega:a^2(R:\omega)=\omega:a^2(\omega:\omega^2)\\
&=\textstyle\frac{1}{a^2}(\omega:(\omega:\omega^2))=\left(\textstyle\frac{\omega}{a}\right)^2
\end{align*}
for all $n>0$.
By substituting $n=1$ and $n=3$, we get 
$\omega^3=a \omega^2$.
By substituting $n=2$, we have $1\in (R:\omega):(R:\omega) = \left(\frac{\omega}{a}\right)^2$. Hence $ a^2\in\omega^2$. 

When this is the case, we further have $R[\frac{\omega}{a}]=R:(R:\omega)\omega=R:\tr_R(\omega_R)$. It follows that 
\begin{align*} 
R:R[\tfrac{\omega}{a}]=&(\tfrac{\omega}{a}:\tfrac{\omega}{a}):R[\tfrac{\omega}{a}]=\tfrac{\omega}{a}:R[\tfrac{\omega}{a}]=\tfrac{\omega}{a}:(R:\tr_R(\omega_R))\\
=&\tfrac{\omega}{a}:(\tfrac{\omega}{a}:\tfrac{\omega}{a}\tr_R(\omega_R))=\tfrac{\omega}{a} \tr_R(\omega_R)=\tfrac{\omega}{a} a (R:\omega)=\tr_R(\omega_R).
\end{align*}
\end{proof}

Theorem \ref{d1.2} is useful to prove Theorem \ref{qqq3.11}.

\section{Almost Gorenstein versus nearly Gorenstein in dimension one}\label{newnewsection2}

Let $(R, \fkm)$ be a Cohen-Macaulay local ring of dimension one. Suppose that $R$ possesses the canonical module $\omega_R$ and a canonical ideal $\omega\subsetneq R$. 
Then, since $\dim R=1$, $\omega$ is an $\fkm$-primary ideal of $R$. Hence we can define the Hilbert function $\ell_R (R/\omega^n)$ of canonical ideal $\omega$. In \cite{CGKM, GMP}, the authors have explored the notions of almost Gorenstein and $2$-almost Gorenstein by using the Hilbert function $\ell_R (R/\omega^n)$. With this background, this section focuses on the case of dimension one. 

First of all, we note the existence of an almost reduction of a canonical ideal (Corollary \ref{jjj2.4}).

\begin{Definition} (\cite[Proposition 1.1]{Lipman} and \cite[before Lemma 8.2]{Mc})\label{g2.1}
Let $R$ be a Noetherian ring and $I$ an ideal of $R$. We then set
\begin{enumerate}[{\rm (a)}] 
\item $\widetilde{I}=\bigcup_{\ell>0} (I^{\ell+1}:_R I^{\ell})$ and
\item $R^I=\bigcup_{\ell>0} (I^{\ell}: I^{\ell})$. 
\end{enumerate}
The ideal $\widetilde{I}$ is called the {\it Ratliff-Rush closure} of $I$, and the ring $R^I$ is coincides with the blow-up of $R$ at $I$ when $\dim R=1$.
\end{Definition}

\begin{Lemma} {\rm (\cite[Lemma 8.2]{Mc})} \label{g2.2}
Let $R$ be a Noetherian ring and $I$ an ideal of $R$ containing a non-zerodivisor of $R$. Then $I^n=(\widetilde{I})^n=\widetilde{I^n}$ for all $n\gg0$.
\end{Lemma}

\begin{Lemma} {\rm (cf. \cite[Proposition 1.1]{Lipman})} \label{jjj2.3}
Let $(R, \fkm)$ be a Cohen-Macaulay local ring of dimension one and $I$ an $\fkm$-primary ideal of $R$. Then we have the following.
\begin{enumerate}[{\rm (a)}] 
\item If there exists a reduction $(a)\subseteq I$ of $I$, then $R^I=R[\frac{I}{a}]=\left(\frac{I}{a}\right)^n$ for all $n\gg 0$.
\item $R^I=I^n:I^n$ for all $n\gg 0$.
\item $R^I=R^{I^n}$ for all $n> 0$.
\item $IR^I\cong R^I$.
\end{enumerate}
\end{Lemma}


\begin{Proposition}\label{g2.5}
Let $(R, \fkm)$ be a Cohen-Macaulay local ring of dimension one and $I$ an $\fkm$-primary ideal of $R$. Then there exist an $\fkm$-primary ideal $J$ and $a\in R$ such that $J\cong I$ and $(a)$ is an almost reduction of $J$.
\end{Proposition}

\begin{proof}
Since Lemma \ref{jjj2.3} (b) and (c), we have $R^I=R^{I^n}=I^n:I^n$ for all $n\gg 0$. Hence $I^n R^{I^n}=(I^n:I^n)I^n=I^n$ is isomorphic to $R^{I^n}=R^I$ by Lemma \ref{jjj2.3} (d). It follows that 
\[
I^{n+1}\cong IR^I\cong R^I \cong I^n.
\]
Hence $I^{n+1}=\alpha I^n$ for some $\alpha \in \rmQ(R)$. Write $\alpha=\frac{a}{b}$ for non-zerodivisors $a$ and $b$ in $R$. Then we obtain $(bI)^{n+1}=a(bI)^n$ and $a\in (bI)^{n+1}:_R (bI)^n\subseteq \widetilde{bI}$. Therefore, we have $a^n \in (\widetilde{bI})^n= (bI)^n$ for $n\gg 0$ by Lemma \ref{g2.2}. 
Thus $(a)$ is an almost reduction of $bI$.
\end{proof}

Due to Proposition \ref{g2.5} we have the following.

\begin{Corollary}\label{jjj2.4}
Let $R$ be a one-dimensional Cohen-Macaulay local ring. If $R$ has a canonical ideal, then $\canred R<\infty$.
\end{Corollary}

While Proposition \ref{g2.5} holds, there exists an $\fkm$-primary ideal which has no parameter almost reduction.

\begin{Example}{\rm (\cite[after Definition 2.1]{Lipman})}\label{kkk2.5}
Let $R=k[[X, Y]]/(XY(X+Y))$, where $k=\mathbb{Z}/2\mathbb{Z}$. Let $x, y$ denote the images of $X, Y$ in $R$. Then any element of $R$ is not an almost reduction of $\fkm=(x, y)$.
\end{Example}

\begin{proof}
Suppose that $\fkm$ has an almost reduction $(a)$. Then, since  $a^n\in \fkm^n \subseteq \fkm$ for $n\gg 0$, $a\in \fkm$. It follows that $(a)$ is a reduction of $\fkm$. 
Write $a=c_1 x+c_2 y+g$, where $c_1, c_2\in k$ and $g\in \fkm^2$. Then we can replace $g$ by $0$. Actually, we have
\[
a\fkm^n\subseteq (a, g) \fkm^n=(a-g)\fkm^n + g \fkm^n\subseteq (a-g)\fkm^n + \fkm^{n+2}\subseteq \fkm^{n+1}
\]
and thus $(a-g)\fkm^n=\fkm^{n+1}$ for $n\gg 0$ by Nakayama's lemma. Hence we may assume that $a$ is either $x$, $y$, or $x+y$. It concludes that $a$ is a zerodivisor of $R$, which is a contradiction since $\fkm^{n+1}\subseteq (a)$.
\end{proof}

Let us continue to explore the canonical reduction number in dimension one.

\begin{Proposition}\label{g2.6}
Let $(R, \fkm)$ be a Cohen-Macaulay local ring of dimension one. Suppose that $R$ possesses a canonical ideal $\omega\subsetneq R$. We then have the following.
\begin{enumerate}[{\rm (a)}] 
\item The following integers are equal:
\begin{enumerate}
\item[{\rm (i)}]  $\canred R$;
\item[{\rm (ii)}]  $\min \left\{ m\ge 0 \mid 
\text{$\ell_R (R/\omega^{n})=\rme_0(\omega)n-\rme_1(\omega)$ holds for all $n\ge m$}\right\}$.
\end{enumerate}
\item $0\le \canred R\le \rme_0(\fkm)-1$, where $\rme_0(\fkm)$ denotes the multiplicity of $R$.
\end{enumerate}
\end{Proposition}

\begin{proof}
(a): 
We may assume that $R$ is not Gorenstein by Proposition \ref{d1.1}. Choose a canonical ideal $\omega \subsetneq R$ so that $\omega$ has  an almost reduction $(a)$. Set $c=\canred R>0$.
Then, for all $n\ge c$, we have 
\begin{align*}
\ell_R(R/\omega^n) &= \ell_R (R/(a)^{n}) - \ell_R (\omega^n/(a)^n) \\
 &= \ell_R (R/(a)){\cdot}n - \ell_R \left(\left(\textstyle\frac{\omega}{a}\right)^n/R\right)\\
 &= \rme_0(\omega)n - \ell_R \left(R\left[\textstyle\frac{\omega}{a}\right]/R\right)
\end{align*}
by Remark \ref{h1.2} (c) and Proposition \ref{f1.2} (a). Note that $R\left[\frac{\omega}{a}\right]$ is independent of the choice of almost reductions by Proposition \ref{f1.2} (b). Thus $c$ is larger than or equals to the integer of (ii). Assume that $c> \text{the integer of (ii)}$. Then 
\begin{align*}
\ell_R(R/\omega^{c-1}) &= \ell_R (R/\omega^c)-\ell_R (\omega^{c-1}/\omega^c)\\
&=\rme_0(\omega)c - \ell_R \left(R\left[\textstyle\frac{\omega}{a}\right]/R\right) - \ell_R \left(\left(\textstyle\frac{\omega}{a}\right)^{c-1}/aR\left[\textstyle\frac{\omega}{a}\right]\right)\\
&>\rme_0(\omega)c - \ell_R \left(R\left[\textstyle\frac{\omega}{a}\right]/R\right) - \ell_R \left(R\left[\textstyle\frac{\omega}{a}\right]/aR\left[\textstyle\frac{\omega}{a}\right]\right).
\end{align*}
On the other hand, we have $\ell_R(R/\omega^{c-1})=\rme_0(\omega)(c-1) - \ell_R \left(R\left[\frac{\omega}{a}\right]/R\right)$. 
It follows that 
\[
\rme_0(\omega)<\ell_R \left(R\left[\textstyle\frac{\omega}{a}\right]/aR\left[\textstyle\frac{\omega}{a}\right]\right)=\rme_0(\omega){\cdot}\rank_R R\left[\textstyle\frac{\omega}{a}\right]=\rme_0(\omega)
\] 
by the multiplicative formula. This is a contradiction.

(b): To prove the inequality we may assume that $R/\fkm$ is infinite by Proposition \ref{i1.7}. Then it follows from (a) and \cite{ES}. 
\end{proof}


Next we study a relation between almost Gorenstein rings and nearly Gorenstein rings in terms of the canonical reduction number. In what follows, throughout this section, let $(R, \fkm)$ be a one-dimensional Cohen-Macaulay local ring with a canonical ideal. Choose a canonical ideal $\omega\subsetneq R$ so that $\omega$ has an almost reduction $(a)$. We start to recall the definitions of almost Gorenstein and nearly Gorenstein.

\begin{Definition}
\begin{enumerate}[{\rm (a)}] 
\item (\cite[Definition 3.1]{GMP}) We say that $R$ is {\it almost Gorenstein} if $\rme_1 (\omega)\le \rmr(R)$. 
\item (\cite[Definition 2.2]{HHS}) We say that $R$ is {\it nearly Gorenstein} if $\tr_R (\omega_R)\supseteq \fkm$.
\end{enumerate}
\end{Definition}

The notion of almost Gorenstein rings originates from Barucci and Fr\"oberg \cite{BF}. They defined almost Gorenstein rings over one-dimensional analytically unramified Cohen-Macaulay rings. After that, Goto, Matsuoka, and Phuong \cite{GMP} defined the notion for arbitrary one-dimensional Cohen-Macaulay local rings, and showed that these definitions are equivalent if there exists a parameter reduction of a canonical ideal, see \cite[Setting 3.4 and Theorem 3.11]{GMP}. Here, we show that we have the equivalence without any assumptions (Proposition \ref{jjj2.7}). Before that, we note the following lemma, which should be known.

\begin{lem}\label{zzzzzz}
Let $I$ and $J$ be fractional ideals of $R$. If $IJ=R$, then $I\cong R$ and $J\cong R$.
\end{lem}

\begin{proof}
Suppose that $IJ=R$, and choose $a_1, \dots, a_n\in I$ and $b_1, \dots, b_n \in J$ such that $\sum_{i=1}^n a_ib_i=1$. Then $a_ib_i \in R$ for all $1 \le i \le n$ since $J\subseteq R:I$ by the assumption $IJ=R$. Hence, by noting that $R$ is local, $a_ib_i$ is a unit for some $1 \le i \le n$. Thus, $a_ib_ic=1$ for some $c\in R$. It follows that $I=(a_i)\cong R$. Indeed, for all $a\in I$, $a=aa_ib_ic=a_i (ab_ic)$ and $ab_ic\in IJ\subseteq I(R:I)\subseteq R$. Hence, $I\subseteq (a_i)$. The reverse inclusion is clear.
\end{proof}

\begin{Proposition}{\rm (generalization of \cite[Theorem 3.11]{GMP})}\label{jjj2.7}
The following conditions are equivalent:
\begin{enumerate}[{\rm (a)}] 
\item $R$ is almost Gorenstein;
\item $\fkm R[\frac{\omega}{a}]= \fkm$;
\item $\fkm \omega \subseteq (a)$.
\end{enumerate}
\end{Proposition}

\begin{proof}
We may assume that $R$ is not Gorenstein. Set $K=\frac{\omega}{a}$. 

(b) $\Leftrightarrow$ (c): It follows from the equivalences
\begin{align*}
\fkm \omega \subseteq (a) & \Leftrightarrow \fkm K \subseteq R  && \Leftrightarrow \fkm K \subsetneq R && \Leftrightarrow \fkm K \subseteq \fkm  \\
& \Leftrightarrow \fkm K^n\subseteq \fkm \text{ for all $n>0$} && \Leftrightarrow \fkm R[K] \subseteq \fkm 
&& \Leftrightarrow \fkm R[K] = \fkm,
\end{align*}
where the second equivalence follows from the fact that $\fkm K=R$ implies $\fkm$ is principal by Lemma \ref{zzzzzz}, that is, $R$ is a discrete valuation ring.

(a) $\Leftrightarrow$ (b): Note that (a) and (b) are equivalent if $R/\fkm$ is infinite, and then $\canred R\le 2$ (\cite[Theorem 3.16 (b)]{GMP}). Hence, by passing to the faithfully flat map $R \to R[X]_{\fkm R[X]}$, we may assume that $\canred R\le 2$ by Proposition \ref{i1.7}.
Note that the proof of Proposition \ref{g2.6} shows that $\rme_1(\omega)=\ell_R(R[\frac{\omega}{a}]/R)=\ell_R(R:\tr_R(\omega_R)/R)$ by Theorem \ref{d1.2}. On the other hand, we have $R:\fkm/R\cong \Ext_R^1(R/\fkm, R)$ by applying the $R$-dual to $0\to \fkm \to R \to R/\fkm \to 0$. Hence we have 
\[
\rme_1(\omega)=\ell_R(R:\tr_R(\omega_R)/R)\ge \ell_R(R:\fkm/R)=\rmr(R)
\] 
since $\tr_R(\omega_R)\subseteq \fkm$. It follows that 
\begin{align*}
& \text{$R$ is almost Gorenstein}  \Leftrightarrow R[\tfrac{\omega}{a}]=R:\tr_R(\omega_R)=R:\fkm \\
\Leftrightarrow & \fkm R[\tfrac{\omega}{a}]\subseteq R \Leftrightarrow \fkm R[\tfrac{\omega}{a}]=\fkm,
\end{align*}
where the third equivalent follows from Lemma \ref{zzzzzz}.
\end{proof}

\begin{Lemma}\label{b1.10}
If $R/(R:R[\frac{\omega}{a}])$ is Gorenstein, then $\canred R\le 2$. 
\end{Lemma}

\begin{proof}
Set $K=\frac{\omega}{a}$ and $S=R[K]$. By applying the $K$-dual $K:-$ to the short exact sequence $0 \to R:S \to R \to R/(R:S) \to 0$, we have 
\[
0\to K \to K:(R:S) \to \Ext_R^1 (R/(R:S), K) \to 0.
\]
On the other hand, we obtain
\[
K:(R:S)=K:((K:K):S)=K:(K:KS)=K:(K:S)=S.
\]
Hence $\omega_{R/(R:S)}\cong \Ext_R^1 (R/(R:S), K) \cong S/K$. 
By our assumption we have $S=~K~+~Rs$ for some $s\in S$. Let $\alpha \in K$. Then, since $sK\subseteq S=K+Rs$, there exists $\beta \in K$ and $r\in R$ such that $s\alpha=\beta+rs$. Whence $s(\alpha-r)=\beta\in K$. It follows that 
\[
(\alpha-r)K^2\subseteq (\alpha-r)S=(\alpha-r)(K+Rs)=(\alpha-r)K+R \beta \subseteq K^2.
\]
Therefore, for $\alpha\in K$ and $\alpha'\in K^2$, there are elements $r\in R$ and $\alpha''\in K^2$ such that $(\alpha-r)\alpha'=\alpha''$. Thus $\alpha \alpha'=r\alpha'+\alpha''\in K^2$. It follows that $K^2=K^3$. Hence, we have $K^2=K^n=R[K]$ for $n\gg 0$, that is, $\canred R\le 2$.
\end{proof}

Now we can illustrate a relation between almost Gorenstein and nearly Gorenstein (see also \cite[Theorem 7.4]{HHS}).

\begin{Theorem}\label{qqq3.11}
The following conditions are equivalent:
\begin{enumerate}[{\rm (a)}] 
\item $R$ is almost Gorenstein;
\item $R$ is nearly Gorenstein and $\canred R\le 2$;
\item $R$ is nearly Gorenstein and $R/(R: R[\frac{\omega}{a}])$ is Gorenstein.
\end{enumerate}
\end{Theorem}

\begin{proof}
For each proof of implication, we may assume that $R$ is not Gorenstein. Set $K=\frac{\omega}{a}$ and $S=R[K]$. 

(a) $\Rightarrow$ (c): By Proposition \ref{jjj2.7}, we have $\fkm S=\fkm\subseteq R$. It follows that $\fkm \subseteq R: S \subseteq R: K \subseteq (R:K)K=\tr_R (\omega_R) \subsetneq R$. Hence $R$ is nearly Gorenstein and $R/(R:S)=R/\fkm$ is Gorenstein.

(c) $\Rightarrow$ (b): It follows from Lemma \ref{b1.10}.

(b) $\Rightarrow$ (a): By Proposition \ref{f1.2} (a), $K^2=S$. It follows that $R:K=(K:K):K=K:K^2=K:S=K:KS=(K:K):S=R:S$.
Hence 
\[
\fkm=\tr_R (\omega_R)=(R:K)K=(R:S)K=R:S,
\]
where the fourth equality follows from $(R:S)K\subseteq (R:S)S=R:S$ and $R:S=(R:S)K^2\subseteq (R:S)K$. Therefore, we get $\fkm K\subseteq \fkm S\subseteq R$.
\end{proof}

Note that the ring in Example \ref{qwqw2.3} is almost Gorenstein which has no parameter reduction (\cite[Example 3.2 (1)]{GMP}). 
In the rest of this section we note an example arising from numerical semigroup rings. 

\begin{Example}
Let $n\ge 3$ and let $R=k[[t^n, t^{n+1}, t^{n^2-n-1}]]$ be a numerical semigroup ring over a field $k$. Then $K=R+Rt$ is a fractional canonical ideal of $R$ such that $R\subseteq K\subseteq \ol{R}$ (\cite[Example (2.1.9)]{GW}). Therefore, we get 
\[
\canred R=\min \left\{ n\ge 0 \mid K^n=K^{n+1} \right\}=n-1.
\] 
On the other hand, one can check that $R$ is a nearly Gorenstein ring for all $n\ge 3$. Hence $R$ is an almost Gorenstein ring if and only if $n=3$.
\end{Example}

\section{The canonical reduction number and idealizations}\label{newsection2}

For a moment, let $R$ be an arbitrary commutative ring and $M$ an $R$-module. 
Let $A=R \ltimes M$ denote the {\it idealization} of $M$ over $R$, that is, $A = R \oplus M$ as an $R$-module and the multiplication in $A$ is given by
$$(a,x)(b,y) = (ab, bx + ay)$$
where $a,b \in R$ and $x,y \in M$. The following are well-known.

\begin{Fact}\label{jjj3.4}
For a local ring $R$ and a nonzero $R$-module $M$, we have the following:
\begin{enumerate}[{\rm (a)}] 
\item (\cite[Theorem 4.8]{AW}) $R\ltimes M$ is a Noetherian ring if and only if $R$ is a Noetherian ring and $M$ is a finitely generated $R$-module.
\item (\cite[Corollary 4.14]{AW}) $R\ltimes M$ is a Cohen-Macaulay ring if and only if $R$ is a Cohen-Macaulay ring and $M$ is a  maximal Cohen-Macaulay $R$-module.
\item (\cite[(7) Theorem]{R}) $R\ltimes M$ is a Gorenstein ring if and only if $R$ is a Cohen-Macaulay ring possessing the canonical module $\omega_R$ and $M\cong \omega_R$.
\end{enumerate}
\end{Fact}

Fact \ref{jjj3.4} shows that the properties of $R\ltimes M$ corresponds to the properties of $R$ and $M$. Especially, if $R$ is a Cohen-Macaulay local ring and $M$ is a maximal Cohen-Macaulay $R$-module, the idealization $R\ltimes M$ builds bridges between  the stratification of Cohen-Macaulay rings and the classification of maximal Cohen-Macaulay modules. 

Set $A=R\ltimes M$. The following theorem gives a characterization of the condition $\canred A\le 2$ via the properties of $M$. Note that we need to assume that $R$ is Gorenstein because of technical reasons, see Propositions \ref{d2.5} and \ref{e2.8}. But, Theorem \ref{d2.3} still provides infinitely many examples of canonical reduction number two in higher dimension. This result also gives better prospects of the result \cite[Corollary 6.4]{GMP}, although there are more general results on the almost Gorenstein property of idealization (\cite[Theorem 6.3]{GMP}, \cite[Theorem 1.1]{GK2}).
Recall that $B$ is a {\it finite birational extension} of $R$ if $B$ is a subring of $\rmQ(R)$ containing $R$ and finitely generated as an $R$-module.

\begin{Theorem}\label{d2.3}
Let $(R, \fkm)$ be a Gorenstein local ring of dimension $d>0$. Let $M$ be a maximal Cohen-Macaulay faithful $R$-module. Set $A=R\ltimes M$. Then the following assertions are equivalent:
\begin{enumerate}[{\rm (a)}] 
\item $\canred A\le 2$;
\item $\tr_A (\omega_A)\cong \omega_A^{*}$;
\item $M$ is isomorphic to some trace ideal of $R$;
\item $M^*$ is isomorphic to some finite birational extension $B$ of $R$ such that $B_\fkM$ is a Cohen-Macaulay local ring of dimension $d$ for all $\fkM\in \Max B$.
\end{enumerate}
When this is the case, if $A$ is not Gorenstein, then $\rmr(A)=\rmr(R/I)+2$ where $I$ denotes the trace ideal isomorphic to $M$.
\end{Theorem}

To prove Theorem \ref{d2.3}, we need some results; the following one should be well-known, but we did not find a good reference.

\begin{Lemma}\label{d2.4}
Let $(R, \fkm)$ be a Noetherian local ring of dimension $d$. Let $B$ be a finite birational extension of $R$. Then the following conditions are equivalent:
\begin{enumerate}[{\rm (a)}] 
\item $B$ is a maximal Cohen-Macaulay $R$-module;
\item $B_\fkM$ is a Cohen-Macaulay local ring of dimension $d$ for all $\fkM \in \Max B$.
\end{enumerate}
\end{Lemma}

\begin{proof}
(a) $\Rightarrow$ (b): Let $\fkM \in \Max B$. Note that $\depth B_\fkM \ge \depth_R B$ since a $B$-regular sequence in $\fkm$ is a $B_\fkM$-regular sequence in $\fkM$.
It follows that 
\[
d=\dim B\ge \dim B_\fkM \ge \depth B_\fkM \ge  \depth_R B=d,
\]
thus $B_\fkM$ is a Cohen-Macaulay local ring of dimension $d$.

(b) $\Rightarrow$ (a): Note that $B$ is a semilocal ring since $B$ is a finite birational extension of $R$. Hence $\fkm B\subseteq J(B)\subseteq \sqrt{\fkm B}$, where $J(B)$ denotes the Jacobson radical of $B$ and $\sqrt{\fkm B}$ denotes the radical of $\fkm B$. Hence 
$\grade(\fkm B, B)=\grade (J(B), B)=d$. It follows that there exists a $B$-regular sequence in $\fkm$ of length $d$.
\end{proof}

The following is a generalization of \cite[Corollary 2.8]{GIK2}. Note that \cite[Corollary 2.8]{GIK2} focuses on $d=1$, but the argument is essentially the same.

\begin{Proposition} \label{d2.5}
Let $(R, \fkm)$ be a Gorenstein local ring of dimension $d>0$. 
Set
\begin{align*}
\varphi: &\left\{ I \ \middle| \   \begin{matrix} 
\text{$I$ is a trace ideal containing a non-zerodivisor of $R$}\\
\text{ and maximal Cohen-Macaulay as an $R$-module}
\end{matrix}\right\}\\
&\to \left\{ B \ \middle| \  \begin{matrix}
\text{$B$ is a finite birational extension of $R$ such that $B_\fkM$ is }\\
\text{a Cohen-Macaulay local ring of dimension $d$ for all $\fkM\in \Max B$}
\end{matrix}
\right\},
\end{align*}	
where $I\mapsto I:I$. Then $\varphi$ is a one-to-one correspondence.
\end{Proposition}

\begin{proof}
(well-definedness): By Remark \ref{b1.7} (b), we have $I:I=R:I\cong I^*$. Hence $I:I$ is a maximal Cohen-Macaulay $R$-module since $R$ is Gorenstein. Hence $(I:I)_\fkM$ is a Cohen-Macaulay ring of dimension $d$ for all $\fkM\in \Max (I:I)$ by Lemma \ref{d2.4}.

(injective): Let $I$ and $J$ be trace ideals containing non-zerodivisors of $R$ and maximal Cohen-Macaulay $R$-modules. If $I:I=J:J$, then we have $I=R:(R:I)=R:(I:I)=R:(J:J)=R:(R:J)=J$.

(surjective): Let $B$ be a finite birational extension of $R$ such that $B_\fkM$ is a Cohen-Macaulay local ring of dimension $d$ for all $\fkM\in \Max B$. Then, by Lemma \ref{d2.4}, $B$ is a maximal Cohen-Macaulay $R$-module. Hence so is $R:B$. $R:B$ is a trace ideal since $\tr_R(B)=(R:B)B=R:B$. Furthermore $(R:B):(R:B)=R:(R:B)B=R:(R:B)=B$ as desired.
\end{proof}



In addition, before proving Theorem \ref{d2.3}, we prove the following proposition.

\begin{Proposition}\label{e2.8}
Let $(R, \fkm)$ be a Gorenstein local ring of dimension $d>1$ and $I$ a trace ideal of $R$. Suppose that $I$ is maximal Cohen-Macaulay as an $R$-module. Then there exists an element $x\in \fkm$ such that $x$ is a non-zerodivisor of $R$ and $R/I$, and $(I+(x))/(x)$ is a trace ideal of $R/(x)$ and maximal Cohen-Macaulay as an $R/(x)$-module.
\end{Proposition}

\begin{proof}
By applying the depth formula to the exact sequence  
\begin{align}\label{ddd4}
0\to I \xrightarrow{\iota} R \to R/I\to 0,
\end{align}
we have $\depth_R R/I\ge d-1>0$. Hence we can choose $x\in \fkm$ so that $x$ is a non-zerodivisor of $R$ and $R/I$. 

Note that the embedding $\iota$ in (\ref{ddd4}) induces the isomorphism 
\begin{align}\label{www4.5}
\Hom_R (I, I)\cong  \Hom_R (I, R)
\end{align}
by \cite[Lemma 2.3]{Lin2} (see also \cite[Proposition 2.1]{GIK2}). Set $\ol{*}=R/(x)\otimes_R *$.  The goal is to prove that the map $\Hom_{\ol{R}} (I \ol{R}, I\ol{R})\to  \Hom_{\ol{R}} (I\ol{R}, \ol{R})$ induced by $\ol{\iota}: I\ol{R} \to \ol{R}$ is bijective. Now we have $\ol{\Hom_R (I, I)}\cong  \ol{\Hom_R (I, R)}$ by (\ref{www4.5}). 
By applying the functor $\Hom_R (I, -)$ to $0 \to R \xrightarrow{x} R \to \ol{R} \to 0$, we obtain $\ol{\Hom_R (I, R)}\cong \Hom_R (I, \ol{R})\cong \Hom_{\ol{R}} (\ol{I}, \ol{R})$. Here, we have $\ol{I}=I/xI=I/((x)\cap I)\cong (I+(x))/(x)=I\ol{R}$. Hence, it is enough to show that we have the natural isomorphism $\ol{\Hom_R (I, I)} \cong \Hom_{\ol{R}} (I \ol{R}, I\ol{R})$. 

By applying the functor $\Hom_R (I, -)$ to $0 \to I \xrightarrow{x} I \to I\ol{R} \to 0$, we get
\begin{align}\label{www4.6}
\begin{split}
0 &\to \Hom_R (I, I) \xrightarrow{x} \Hom_R (I, I) \to \Hom_R (I, I\ol{R}) \\
& \to \Ext_R^1 (I, I) \xrightarrow{x} \Ext_R^1 (I, I).
\end{split}
\end{align}
On the other hand, by applying the functor $\Hom_R (I, -)$ to (\ref{ddd4}), we have $\Hom_R (I, R/I)\cong \Ext_R^1 (I, I)$ by (\ref{www4.5}). It follows that $x$ is a non-zerodivisor of $\Ext_R^1 (I, I)$ since 
\[
\Ass_R (\Ext_R^1 (I, I))=\Supp_R I \cap \Ass_R (R/I) \subseteq \Ass_R (R/I).
\] 
Hence (\ref{www4.6}) provides the isomorphism $\ol{\Hom_R (I, I)} \cong \Hom_R (I, I\ol{R}) \cong \Hom_{\ol{R}} (I \ol{R}, I\ol{R})$ as desired.
\end{proof}

Now let us prove Theorem \ref{d2.3}.

\begin{proof}[Proof of Theorem \ref{d2.3}]
(a) $\Leftrightarrow$ (b) follows from Theorem \ref{d1.2}, and (c) $\Leftrightarrow$ (d) follows from Proposition \ref{d2.5}.
Thus we have only to show that (b) $\Leftrightarrow$ (c).

(b) $\Rightarrow$ (c): By Corollary \ref{d1.6}, we may assume that $A$ is generically Gorenstein. 
For each $\fkp\in \Min R$, we have $P=\fkp\times M \in \Min A$, whence $A_P\cong R_\fkp\ltimes M_\fkp$ is Gorenstein. Hence, by Fact \ref{jjj3.4} (c), $M_\fkp\cong R_\fkp$ since $M$ is faithful. It follows that $M$ is of rank one. Since $M\cong M^{**}$ is torsionfree, $M\cong I$ for some ideal $I$ of $R$.
We may assume that $M=I$. Then $A=R\ltimes I$. Set 
\[
K=(R:I) \times R
\]
as an $R$-module. For $(a, x)\in A$ and $(b, y)\in K$, let us define an $A$-action into $K$ as follows:
\[
(a, x){\circ}(b, y)=(ab, ay+bx).
\]
With this action, $K$ is an $A$-module. It is standard to show that 
\[
K\cong \Hom_R(A, R) \cong \omega_A
\]
as $A$-modules. Furthermore we have $A\subseteq K\subseteq \rmQ(A)=\rmQ(R)\times \rmQ(R)$, thus $K$ is a fractional canonical ideal of $A$.

On the one hand, we get
\begin{align*}
\omega_A^* \cong A:K=(R\ltimes I): ((R:I)\times R)=I\times (I:(R:I)),
\end{align*}
where the last equality follows from the following argument. Let $(a, x)\in \rmQ(A)=\rmQ(R)\times \rmQ(R)$. Then
\begin{align*}
&(a, x) \in (R\ltimes I): ((R:I)\times R) \\
 \Leftrightarrow& (a,x)(b,y)=(ab,ay+bx)\in R\ltimes I \quad \text{for all $b\in R:I$ and all $y\in R$}\\
 \Leftrightarrow& \begin{cases}
a\in R:(R:I)=I,\\
x\in I:(R:I).
\end{cases}
\end{align*}

On the other hand, we get
\begin{align*}
\tr_A(\omega_A)=&(A:K)K\\
=&(I\times (I:(R:I))) ((R:I)\times R)\\
=&(R:I)I\times (I+(R:I)(I:(R:I)))\\
=& (R:I)I \times I.
\end{align*}

Hence $\tr_A (\omega_A)\cong \omega_A^{*}$ shows that $I\cong (R:I)I$. It follows that $M=I\cong \tr_R(I)$.

(c) $\Rightarrow$ (b): Let $M\cong I$ for some trace ideal $I$ of $R$. Then $M_\fkp\cong I_\fkp=R_\fkp$ for all $\fkp\in \Min R$ since $M$ is faithful. By noting that $\Min A=\{\fkp\times M \mid \fkp\in \Min R\}$ since $(0\times M)^2=0$, $A$ is generically Gorenstein. We may assume that $M=I$. Then, by the same argument of the proof of (a) $\Rightarrow$ (c), $K=(R:I) \times R$ is a fractional canonical ideal of $A$. Hence, to prove $\tr_A (\omega_A)\cong \omega_A^*$, it is enough to show $I=(R:I)I$ and $I:(R:I)=I$ by the proof of (a) $\Rightarrow$ (c).

We obtain that $(R:I)I=I$ by Remark \ref{b1.7} (b). Furthermore, by Remark \ref{b1.7} (b), we have $I\subseteq I:(I:I)=I:(R:I)\subseteq R:(R:I)=I$.
Therefore, we have 
\[
\tr_A(\omega_A)= (R:I)I \times I=I\times (I:(R:I)) \cong \omega_A^* 
\]
as desired.

When this is the case, assume that $A$ is not Gorenstein. Let $I$ be the trace ideal isomorphic to $M$. In the calculation of $\rmr(A)$, we may assume that $M=I$.
Choose a canonical ideal $\omega \subsetneq A$ and  $a \in A$ so that $\omega^3=a\omega^2$ and $a^2\in \omega^2$. Then $A/\omega$ is Gorenstein of dimension $d-1$ (Fact \ref{f1.5}). Hence we can choose $x_1, \dots, x_{d-1}\in \fkm$ so that $x_1, \dots, x_{d-1}$ is an $R$, $R/I$, $A$, and $A/\omega$-sequence. 
Then, by  Proposition \ref{e2.8}, we can pass to $R\to R/(x_1, \dots, x_{d-1})$. Thus we may assume that $d=1$. In the case of dimension one, our assertion follows from  \cite[Proposition 6.5]{GKL}.
\end{proof}

\begin{Corollary}\label{d2.6}
Let $(R, \fkm)$ be a Gorenstein local domain of dimension $d>0$ and $M$ a nonzero maximal Cohen-Macaulay $R$-module. Set $A=R\ltimes M$. Then the following assertions are equivalent:
\begin{enumerate}[{\rm (a)}] 
\item $\canred A\le 2$;
\item $\tr_A (\omega_A)\cong \omega_A^{*}$;
\item $M$ is isomorphic to some trace ideal of $R$;
\item $M^*$ is isomorphic to some finite birational extension $B$ of $R$ such that $B_\fkM$ is a Cohen-Macaulay local ring of dimension $d$ for all $\fkM\in \Max B$.
\end{enumerate}
When this is the case, if $A$ is not Gorenstein, then $\rmr(A)=\rmr(R/I)+2$ where $I$ denotes the trace ideal isomorphic to $M$.
\end{Corollary}

\begin{proof}
Let $M$ be a nonzero maximal Cohen-Macaulay $R$-module. Then, since $M\cong M^{**}$, $M$ can be embedded into a finitely generated free $R$-module $F$. Hence, if $aM=0$ for $a\in R$, then $a=0$ since $R$ is a domain.
\end{proof}

We close this paper with examples of Corollary \ref{d2.6} arising from semigroup rings.

\begin{Definition}
Let $\mathbf{a}_1, \mathbf{a}_2, \ldots, \mathbf{a}_\ell \in \Bbb Z^n~(\ell >0)$ be lattice points. Set 
\[
C = \left<\mathbf{a}_1, \mathbf{a}_2, \ldots, \mathbf{a}_\ell \right>=\left\{\sum_{i=1}^\ell c_i\mathbf{a}_i \ \middle| \  0 \le c_i \in \Bbb Z~\text{for~all}~1 \le i \le \ell \right\}
\]
and call it the {\it semigroup} generated by $\mathbf{a}_1, \mathbf{a}_2, \ldots, \mathbf{a}_\ell$. Let $S=k[[X_1, \dots, X_n]]$ be the formal power series ring over a field $k$. We then set
\[
k[[C]] = k[[\mathbf{X}^{\mathbf{a}_1}, \mathbf{X}^{\mathbf{a}_2}, \ldots, \mathbf{X}^{\mathbf{a}_\ell}]]
\]
in $S$, where $\mathbf{X}^{\mathbf{a}}=X_1^{a_1} X_2^{a_2} \cdots X_n^{a_n}$ for $\mathbf{a}=(a_1, a_2, \dots, a_n)$.
The ring $k[[C]]$ is called the {\it semigroup ring} of $C$ over $k$. 
\end{Definition}

\begin{Proposition}
Let $\mathbf{a}_1, \mathbf{a}_2, \ldots, \mathbf{a}_\ell \in \Bbb Z^n~(\ell >0)$ be lattice points.
For a positive integer $m$, set 
\[
R_m=k[[\mathbf{X}^{m\mathbf{a}_1}, \mathbf{X}^{m\mathbf{a}_2}, \ldots, \mathbf{X}^{m\mathbf{a}_\ell}]].
\]
Suppose that $R_1$ is a Gorenstein normal domain of dimension $n$. Then $R_m$ is a Gorenstein domain and the integral closure of $R_m$ is $R_1$. Furthermore, if $m=ab$ with some positive integers $a$ and $b$, then 
\[
R_m\subseteq R_a \subseteq R_1
\] 
and $R_a$ is a finitely generated  $R_m$-module. 
Therefore, the canonical reduction number of $R_m\ltimes \Hom_{R_m} (R_m, R_a)$ is two.
\end{Proposition}

\begin{proof}
Note that the $k$-algebra homomorphism $k[[X_1, \dots, X_n]] \to k[[X_1, \dots, X_n]]$, where $X_i \mapsto X_i^m$ for $1\le i \le n$, provides the isomorphism $R_1\cong R_m$ of rings. Thus $R_m$ is Gorenstein. 
Since $(\mathbf{X}^{\mathbf{a_j}})^m\in R_m$ for all $1\le j \le \ell$, we get 
\[
R_m\subseteq R_1 \subseteq \ol{R_m}\subseteq \ol{R_1},
\] 
where $\ol{R}$ denotes the integral closure of $R$. Thus $\ol{R_m}=R_1$ since $R_1$ is normal. 

Suppose that $m=ab$. Then we have $R_m\subseteq R_a \subseteq R_1$ and $R_a$ is finitely generated as an $R_m$-module by the following claim.
\end{proof}

\begin{claim}\label{claim1}
\[
R_a=\sum_{\begin{smallmatrix}
\mathbf{c}=c_1 (a \mathbf{a}_1) +  \cdots + c_\ell (a \mathbf{a}_\ell)\\
\text{for $0\le c_1, c_2, \dots, c_\ell<b$}
\end{smallmatrix}} R_m \mathbf{X}^{\mathbf{c}}
\]
\end{claim}

\begin{proof}[Proof of Claim \ref{claim1}]
The inclusion $\supseteq$ is clear, thus we have only to show the inclusion $\subseteq$. Let $\mathbf{X}^{\mathbf{c}}\in R_a$. Write 
$\mathbf{c}=c_1 (a \mathbf{a}_1) +  \cdots + c_\ell (a \mathbf{a}_\ell)$, where $0\le c_1, c_2, \dots, c_\ell$. Then we can find integers $r_j$ and $0\le q_j<b$ such that $c_j=br_j+q_j$ for all $1\le j\le \ell$. Hence we get
\begin{align*}
\mathbf{c}=\sum_{j=1}^\ell c_j a \mathbf{a}_j=\sum_{j=1}^\ell (br_j+q_j) a \mathbf{a}_j=\sum_{j=1}^\ell (r_j m \mathbf{a}_j +  q_j a \mathbf{a}_j).
\end{align*}
It follows that $\mathbf{X}^{\mathbf{c}}= \left(\mathbf{X}^{\sum_{j=1}^\ell r_j m\mathbf{a}_j}\right) \left(\mathbf{X}^{\sum_{j=1}^\ell q_j a \mathbf{a}_j}\right)$ is in the right hand side of Claim \ref{claim1}.
\end{proof}

\begin{Example}
For $m>0$, let $R_m=k[[X^m, X^mY^m, X^mY^{2m}]]$ be a semigroup ring over a field $k$. Then $\canred (R_m\ltimes \Hom_{R_m} (R_m, R_a))=2$ if $a$ divides $m$, since $R_1$ is a Gorenstein normal domain by \cite[Theorem 6.3.5]{BH}.
\end{Example}

\begin{acknowledgments}
The author thanks Ryotaro Isobe for telling him Example \ref{kkk2.5}. The author also thank the anonymous referee for reading the paper carefully and giving him helpful comments.
\end{acknowledgments}



\end{document}